\documentclass[12pt, 14paper,reqno]{amsart}
\usepackage{amsmath,amsfonts,amssymb}

\theoremstyle{plain}
\usepackage{color}
\newtheorem{theorem}{Theorem}

\theoremstyle{proof}
\theoremstyle{definition}

\theoremstyle{remark}

\theoremstyle{lamma}

\numberwithin{equation}{section}
\numberwithin{lemma}{section}
\numberwithin{theorem}{section}

\usepackage{amsmath}
\usepackage{amsfonts}   
\usepackage{amssymb}
\usepackage{amssymb, amsmath, amsthm}
\usepackage[breaklinks]{hyperref}
\theoremstyle{thmrm}

\usepackage{graphicx}
\begin{document} 
\title[On the solutions of  a Lebesgue - Nagell type  equation]{On the solutions of  a Lebesgue - Nagell type  equation}
\author{Sanjay Bhatter, Azizul Hoque and Richa Sharma}
\address{Sanjay Bhatter @Sanjay Bhatter, Malaviya National Institute of Technology Jaipur,  Jaipur, 302017, India.}
\email{bhatters@gmail.com}
\address{Azizul Hoque @Azizul Hoque Harish-Chandra Research Institute,
Chhatnag Road, Jhunsi,  Allahabad 211 019, India.}
\email{ azizulhoque@hri.res.in}
\address{Richa Sharma @Richa Sharma, Malaviya National Institute of Technology Jaipur,  Jaipur, 302017, India.}
\email{richasharma582@gmail.com}

\maketitle

\begin{abstract}
We find all positive integer solutions in $x, y$ and $n$ of $x^2+19^{2k+1}=4y^{n}$ for any non-negative integer $k$.
\end{abstract}

\section{Introduction}
Lebesgue - Nagell type equations are of the form
\begin{equation}\label{LNTE}
x^2+D=\lambda y^n,
\end{equation}
where $D$ and $\lambda$ are fixed positive integers, and a solution is given by a triple $(x, y, n)$ of positive integers. Many special cases of \eqref{LNTE} have been considered over the years, but all most all the results for general values of $n$ are of fairly recent origin. The earliest result is due to P. de Fermat who had shown that for $n=3$ \eqref{LNTE} has only one solution, viz. $(x,y)=(5,3)$ 
when $\lambda=1$ and $D=2$. There are many interesting results on the integer solutions of \eqref{LNTE} when $\lambda=1$. The first result in this case is due to V. A. Lebesgue \cite{LE1850}, who proved that \eqref{LNTE} has no solutions in positive integers $x, y$ and $n$ when $D=1$. For $D > 1$ and for $\lambda=1$, \eqref{LNTE} has been extensively studied by several authors and in particular, by J. H. E. Cohn \cite{JC03, JC93}, A. Hoque and H. K. Saikia \cite{HS16}, and M. Le \cite{LM93, LM95}. We also refer to \cite{JC93} for a complete survey on \eqref{LNTE} when $\lambda=1$. In \cite{JC93}, J. H. E. Cohn solved \eqref{LNTE} completely for $77$ values of $D \leq 100$ in the case of $\lambda=1$. M. Mignotte and B. M. M. de Weger \cite{MW96} solved \eqref{LNTE} for $D = 74, 86$ and $\lambda=1$. Also M. A. Bennett and C. M. Skinner \cite{BS04} treated \eqref{LNTE} when $D = 55, 95$ and $\lambda=1$.  
Recently, Y. Bugeaud, M. Mignotte and S. Siksek \cite{BM20} studied \eqref{LNTE} completely for the remaining $19$ values of $D\leq 100$ when $\lambda=1$.
 
The equation \eqref{LNTE} is naturally well connected with the investigation of the class number of the imaginary quadratic number field $\mathbb{Q}(\sqrt{-D})$. The solvability of some special cases of \eqref{LNTE} has been used in \cite{CHKP, CH} to investigate the class numbers of certain imaginary quadratic number fields. In \cite{AM}, S. A. Arif and F. S. A. Muriefah proved for the case $(D, \lambda)=(3^{2\ell+1}, 1)$ that \eqref{LNTE} has only one solution, given by $(x,y,n)=(10\times 33^m, 7\times 3^{2m}, 3)$ only when $\ell=2+3m$ with some integer $m\geq 0$. F. Luca \cite{LU} extended this result, and solved \eqref{LNTE} completely for the case $(D,\lambda)=(3^\ell, 1)$. This result was further generalized by S. A. Arif and F. S. A. Muriefah in \cite{AM02}, and A. B\'{e}rczes and I. Pink in \cite{BP08} to the case $(D, \lambda)=(p^\ell, 1)$ for certain primes $p$ and integers $\ell$. In \cite{SS, HM11}, the authors considered \eqref{LNTE} for the case $(D,\lambda)=(d^{2\ell+1}, 1)$, where $d>0$ is a square-free integer, and solved it completely under the assumption that the class number of $\mathbb{Q}(\sqrt{-d})$ is $1$ apart from the case $d \equiv 7 \pmod 8$ in which $y$ was supposed to be odd.  In a recent work, A. B\'{e}rczes and I. Pink \cite{BP} extended the result of \cite{SS} to the case: class number of $\mathbb{Q}(\sqrt{-d})$ is either $2$ or $3$.  
For the case $\lambda=2$, Sz. Tengely \cite{TE} considered \eqref{LNTE} when $D=p^{2q}$ with $p$ and $q>3$ odd primes, and proved that it has finitely many solutions in $x, y, p, q$ under the assumption that $y$ is restricted to the set of integers which are not the sum of two consecutive squares. A more general version of this case was closely studied by F. S. A. Muriefah, F. Luca, S. Siksek and Sz. Tengely in \cite{MFST} when $D \equiv 1 \pmod 4$. On the other hand for the case $\lambda=4$, F. Luca, Sz. Tengely and A. Togb\'{e} \cite{LTT} studied the solutions in positive integers $x, y, n$ of \eqref{LNTE} when $D\equiv 3\pmod 4$ with $D\leq 100$. They also determined all possible solutions in positive integers $x, y, n$ when $D = 7^a\times 11^b\times 13^c$, with $a, b, c \geq 0$ such that $\min\{b, c\} = 0$.

In this paper, we investigate the solutions $(x,y,n)$ of \eqref{LNTE} in positive integers  
when $(D, \lambda) =(19^{2k+1},4)$ for any non-negative integer $k$.
More precisely, we prove:

\begin{theorem}\label{thm} Let $k\geq 0$ be an integer. 
The Diophantine equation
\begin{equation}\label{1}
x^2+19^{2k+1}=4y^{n}
\end{equation}
has no solutions in positive integers $x, y, n$ except $$(x,y,n)\in\left\{\left(19^t\times \frac{19^{2(k-t)+1}-1}{2},19^t\times\frac{19^{2(k-t)+1}+1}{4}, 2\right), (559\times 19^{7m},5\times 19^{2m},7)\right\}$$ with $t, m\in \mathbb{Z}_{\geq 0}$ satisfying $k=7m$ in case of the second solution, and $n \neq 1 $. For $n=1$, it has infinitely many solutions in positive integers $x,y, n$.
\end{theorem}

We use elementary arguments and some properties of Lucas numbers. Let us suppose $\alpha$ and $\beta$ be two algebraic numbers satisfying:
\begin{itemize}
\item $\alpha +\beta$ and $\alpha \beta$ are nonzero coprime rational integers,
\item $\frac{\alpha}{\beta}$ is not a root of unity.
\end{itemize}
 Then the sequence $(u_{n})_{n=0}^{\infty}$ defined by
 $$
 u_{n} = \frac{\alpha^{n}- \beta^{n}}{\alpha-\beta}, \quad n\geq 1,
 $$ 
is a Lucas sequence. If, instead of supposing that $ \alpha + \beta\in \mathbb{Z}$ , we only suppose that $(\alpha+\beta)^{2}$ is a non-zero rational integer coprime prime to $\alpha \beta$, then the Lehmer sequence $(u_{n})_{n=0}^{\infty}$ associated to $\alpha$ and $\beta$ is defined by
$$
u_{n}=\begin{cases}
\frac{\alpha^{n}-\beta^{n}}{\alpha-\beta} \quad
\text{if } n \text{ is odd} , \\ 
 \frac{\alpha^{n}-\beta^{n}}{\alpha^{2}-\beta^{2}} \quad
 \text{if } n \text{ is even}. 
  \end{cases} 
 $$
 We say that a prime number $p$ is a primitive divisor of a Lucas number $u_{n}$ if $p$ divides $u_{n}$ but does not divide $(\alpha-\beta)^{2}u_{2} \cdots u_{n-1}.$ Similarly, $p$ is a primitive divisor of a Lehmer number $u_{n}$ if $p$ divides $u_{n}$ but not $(\alpha^{2}-\beta^{2})^{2}u_{3}\cdots u_{n-1}$.
\section{Proof of Theorem \ref{thm}}
It is clear from \eqref{1} that both $x$ and $y$ are odd for $n\geq 1$. We first treat the case when $n=1$. In this case, \eqref{1} becomes 
\begin{equation}\label{eq1}
x^2+19^{2k+1}=4y.
\end{equation}
As $x$ is odd, let $x=2t+1$ for any non-negative integer $t$. Then \eqref{eq1} gives $y=t^2+t+\frac{1+19^{2k+1}}{4}$. Since $\frac{1+19^{2k+1}}{4}$ is an integer, so that  \eqref{eq1} has infinitely many solutions in positive integers $x$ and $y$ due to the infinitely many choice of $t$. In fact, all these solutions are given by the following parametric form: 
\begin{eqnarray*}
\begin{cases}
x=2t+1, \\ y=t^2+t+\frac{1+19^{2k+1}}{4}, \text{ where}~ \ t \in  \mathbb{Z}_{\geq 0}.
\end{cases}
\end{eqnarray*}
Next we look at the case when $n$ is even and $19\nmid x$. In this case, \eqref{1} can be written as 
\begin{equation}\label{x}
19^{2k+1}=(2y^{m}-x)(2y^{m}+x),
\end{equation}
where $n=2m$. If $\gcd (2y^m-x, 2y^m+x)\ne 1$, then \eqref{x} gives that $19$ divides both $2y^m-x$ and $2y^m+x$ which imply $19\mid x$. This is a contradiction. Thus $\gcd (2y^m-x, 2y^m+x)=1$, and hence \eqref{x} gives 
\begin{equation}\label{a}
2y^m-x=1
\end{equation}
and 
\begin{equation}\label{aa}
2y^m+x=19^{2k+1}.
\end{equation}
Substrating \eqref{aa} from \eqref{a}, and then reading modulo $3$, we see that $x\equiv 0\pmod 3$. Further reading \eqref{a} modulo $3$, we see that $y^m\equiv 2\pmod 3$. This shows that $m$ is odd and $y\equiv 2\pmod 3$. Thus \eqref{x} can be written as 
$$x^2+19^{2k+1}=4Y^m$$ 
with $Y=y^2$. This is treated by the next case, and we see that it does not contribute any solution in positive integers $x, y, n$, except the case $m=1$. For $m=1$, \eqref{a} and \eqref{aa} together imply 
$$(x, y)=\left(\frac{19^{2k+1}-1}{2},\frac{19^{2k+1}+1}{4}\right).$$

We now treat the case when $n$ is odd and $19\nmid x$. For this case, we can consider without loss of generality that $n=p$ for some odd prime $p$. Thus \eqref{1} becomes,
\begin{align}\label{xxx}
 x^{2} + 19^{2k+1} = 4y^{p}.
\end{align}
As the class number of $\mathbb{Q} (\sqrt{-19}) $ is $1$, we have the following factorization,
$$
(x+ 19^{k}\sqrt{-19})(x- 19^{k}\sqrt{-19})=4y^{p}.
$$
This can be rewritten as:
$$
\bigg(\frac{x+ 19^{k}\sqrt{-19}}{2}\bigg) \bigg(\frac{x- 19^{k}\sqrt{-19}}{2}\bigg)=y^p.
$$
It is easy to see that $\gcd \bigg(\frac{x+ 19^{k}\sqrt{-19}}{2}, \frac{x- 19^{k}\sqrt{-19}}{2}\bigg)=1$ since $19\nmid x$. Thus we can write 
\begin{equation}\label{eq2}
\left(\frac{x+ 19^{k}\sqrt{-19}}{2}\right)= \left(\frac{a+b\sqrt{-19}}{2}\right)^p ,\quad a\equiv b\pmod 2
\end{equation}
satisfying
$$
y=\frac{a^{2}+19b^{2}}{4}.
$$
Since the only units in the ring of integers of $\mathbb{Q}(\sqrt{-19})$ are $\pm 1$, these can be absorbed into the $p$-th power.

If both $a$ and $b$ are even then $y$ is also even which is a contradiction. Thus both $a$ and $b$ are odd integers. 

Equating the imaginary parts in \eqref{eq2}, we obtain
\begin{align} \label{eq3}
2^{p-1} (19)^{k}=b \sum_{r=0}^{\frac{p-1}{2}} \binom{p}{2r+1}a^{p-2r-1}(-19)^{r}(b^{2})^{r}
\end{align}
with $\gcd(19,a)=1$.
We now discuss \eqref{eq3} in three cases.\\
Case i. $b=\pm 1$. In this case \eqref{eq3} gives
\begin{equation}\label{xy}
\pm 2^{p-1} (19)^{k}= \sum_{r=0}^{\frac{p-1}{2}} \binom{p}{2r+1}a^{p-2r-1}(-19)^{r}.
\end{equation}
If $k=0$, then \eqref{xy} becomes
\begin{align}  \label{12}
\pm 2^{p-1} = \sum_{r=0}^{\frac{p-1}{2}} \binom{p}{2r+1}a^{p-2r-1}(-19)^{r}.
\end{align}
We will treat \eqref{12} later.

Let us suppose $k>0$. Then it is easy to see from \eqref{xy} that $p=19$ since $\gcd(19,a)=1$. Thus dividing both sides of \eqref{xy} by $19$, and then reading modulo $19$, we deduce that $k=1$ and only the positive sign holds. More precisely, we have 
\begin{align} \label{eq4}
 2^{p-1}=a^{p-1}-\sum_{r=1}^{\frac{p-1}{2}} \binom{p}{2r+1} a^{p-2r-1}(-19)^{r-1}.
 \end{align}
 This is a particular case of \eqref{eq5}.
 \\
Case ii. $b=\pm 19^{t}$ with $0<t<k$. Then \eqref{eq3} implies  
$$
\pm 2^{p-1} 19^{k-t}= \sum_{r=1}^{\frac{p-1}{2}} \binom{p}{2r+1} a^{p-2r-1}(-19)^{r}(19^{2t})^{r}.
$$ 
We reduce this equation modulo $19$ and derive as before that,
$p=19$. Then dividing by $19$, we note that $k-t-1=0 $ and only the positive sign holds. Consequently, 
\begin{align} \label{eq5}
 2^{p-1}=a^{p-1}-\sum_{r=1}^{\frac{p-1}{2}} \binom{p}{2r+1} a^{p-2r-1}(-19)^{r-1}(19^{2r})^{k-1}.
 \end{align}
 It is noted when $k=1$, \eqref{eq4} is same as \eqref{eq5}.
 
 We now prove that \eqref{eq5} does not hold when $p=19$.
 Here, $p\equiv 3 \pmod 4$ so that $p= 3+2^s m $, where $s\geq 2$ and $(2,m)=1$. \\ 
 Considering \eqref{eq5} modulo $2^{s+1}$,
 we get 
 $$
 0\equiv a^{p-1}-\binom{p}{3}a^{p-3}19^{2t}+(-19)^{\frac{p-1}{2}}(19^{2t})^{\frac{p-1}{2}}(19)^{-1} \pmod {2^{s+1}}.
 $$
 The first term is (using Euler's theorem)
 $$
 a^{p-1}\equiv a^{2} \pmod {2^{s+1}}.
 $$
 Similarly, the second term is  
 $$
 \binom{p}{3} a^{p-3}\equiv 1+2^{s-1}m \pmod {2^{s+1}}\equiv 1+2^{s-1} \pmod {2^{s+1}}.
 $$
 Using this as before with the  observation that $ 2^{a}t\equiv 2^{a}\pmod {2^{a+1}}$  if $t$ is odd, 
 the third term is
 $$ 
 (-19)^{\frac{p-1}{2}}(19^{2t})^{\frac{p-1}{2}}(19)^{-1}\equiv -19^{2t} \pmod {2^{s+1}},~ \text{where}~ t> 0.
 $$
 Hence 
  $$
  0\equiv a^{2}-(1+2^{s-1})19^{2t}-19^{2t} \pmod {2^{s+1}}.
  $$
 This implies that 
 $$
 a^{2} \equiv 19^{2t}(2+2^{s-1}) \pmod {2^{s+1}}
 $$ 
 which is not possible since $a$ is odd and $\gcd(19,a)=1$.

Case iii. $b=\pm 19^{k}$. In this case \eqref{eq3} gives
\begin{align} \label{3}
\pm 2^{p-1} = \sum_{r=0}^{\frac{p-1}{2}} \binom{p}{2r+1} a^{p-2r-1}(-19)^{r}(19^{2k})^{r}.
\end{align}
If $k=0$, we get \eqref{12}. 
Earlier, we have shown that to solve \eqref{1}, it suffices to consider \eqref{eq3}, where $a$ is odd integer, $b=\pm 19^{k}$ and $a$, $p$ and $19$ satisfy \eqref{3}.\\

Let $ \alpha =\frac{a+ 19^{k} \sqrt{-19}}{2}$. Then  
 \begin{align} \label{4}
\frac{\alpha^{p}-\bar{\alpha}^{p}}{\alpha-\bar{\alpha}}=\pm 1,
\end{align}
where $\bar{\alpha}$ is the conjugate of $\alpha$. It is obvious that $\alpha$ and $\bar{\alpha}$ are algebraic integers. We also have $\alpha+ \bar{\alpha}= a$ and  $\alpha \bar{\alpha}=\frac{a^{2}+19^{2k+1}}{4}=y$ since $b=\pm 19^{k}$. Furthermore , since $\frac{a^{2}+19^{2k+1}}{4}=y$, so that $19^{2k+1}\mid \gcd(y,a)$. Thus $19\nmid y$ and hence $\gcd(y,a)=1.$ 
Therefore $\alpha + \bar{\alpha}$ and $\alpha \bar{\alpha}$ are coprime. We also observe that $\frac{\alpha}{\bar{\alpha}} \not \in \{\pm 1\}$, the set of roots of unity in $\mathbb{Q} (\sqrt{-19}) $. Thus $(\alpha , \bar{\alpha})$ is a Lucas pair and hence
$u_{p}=\frac{\alpha^{p}-\bar{\alpha}^{p}}{\alpha-\bar{\alpha}}$ is a Lucas number. Using \cite{BH01}, we see that $u_{p}$ has primitive divisors for all primes $p>13$, which contradicts to \eqref{4} and thus \eqref{eq3} has no solution for all primes $p>13$.

Also if $ p \in  \{5,7,11,13 \}$, then there are Lucas pairs $(\alpha,\bar{\alpha})$ for which $u_{p}$ does not have primitive divisors.
We consider each $p$ separately.

For $p =13$, the only Lucas pair without primitive divisors  belongs to $\mathbb{Q} (\sqrt{-7}) $ and so it is not possible.

Also for $ p =11$, there is no Lucas pair without primitive divisors. Thus  there do not exist any solution for $p = 11, 13$.

Again for $p=7$, the only Lucas pair of the required form is $(\alpha, \bar{\alpha})$ where $ \alpha = \frac{1-\sqrt{-19}}{2}.$
In such a situation we have $a=1$ and $k=0$. Now using \eqref{eq3}, we see that the only possible solution is $(x,y,p)=(559,5,7).$ Thus we get a solution of \eqref{xxx} only when $k=0$ which is  $(x,y,p)=(559,5,7)$.
 
 Further for $p=5$, there is no Lucas pair $({\alpha,\bar{\alpha}})$ without primitive divisors of the desired form
 $ \alpha= a + 19^{k} \sqrt{-19}$ with odd $a$. 
 
Finally for $p=3$, equating real and imaginary parts in \eqref{eq2}, we obtain
\begin{equation}\label{p31}
4x=a^3-57ab^2
\end{equation}
and 
\begin{equation}\label{p32}
4\times 19^k=3a^b-19b^3.
\end{equation}
If both $a$ and $b$ are even then \eqref{p31} implies $2\mid x$ which is a contradiction. Again if both $a$ and $b$ are odd then reading \eqref{p32} modulo $3$, we obtain $b\equiv 2  \pmod 3$. Thus we can write $b=3r+2$ for some integer $r$, and then replacing this value in \eqref{p32}, we get
$$4\times 19^k=9a^2r+6a^2-19(27r^3+54r^2+36r+8).$$  
Reading this modulo $9$, we arrive at
$4\equiv 6a^2-8\pmod 9$. This further implies $a^2\equiv 2\pmod 3$ which is not possible.
 Thus we complete all the above cases.
 
 The only remaining case to treat is $19\mid x$. 
Let $x=19^{s}X$, and $y=19^{t}Y$ such that $\gcd(19,X)=\gcd(19,Y)=1$ with positive integers $ s$ and $t$. Then \eqref{1} becomes
 \begin{equation}\label{xx}
  19^{2s}X^{2}+19^{2k+1}=4(19)^{tn}Y^{n}.
 \end{equation}
  We encounter this equation with three possibilities.
 Let $2k+1= \mbox{min} \{2s,2k+1,t n\}$. Then \eqref{xx}
  $$
  19(19^{s-k-1}X)^{2}+1=4Y^{n}(19)^{tn-2k-1}.
  $$
 Again utilising previous technique we read it modulo $19$, and conclude that $t n - 2 k -1= 0$. Thus
 $$
 19(19^{s-k-1}X)^{2}+1=4Y^{n}.
 $$
This equation has no solution using \cite{LM93}.
  
  If
 $nt= \mbox{min} \{2s,2k+1,n t\}$. Then
 $$
 19^{2s - n t}X^{2}+19^{2k- n t+1}=4Y^{n}.
 $$
 This relation holds only if $2 k - n  t+1=0$ or if $2 s=n t$.
 If $2k-nt+1=0$, then 
 $$
 19^{2 s - n t}X^{2}+1=4Y^{n}.
 $$
 Thus as in the previous case, it has no solution.
  Again if $2 s = n t$, then
 $$
 X^{2}+19^{2(k-s)+1}=4Y^{n}.
 $$
 This equation has no solution by the previous cases, except the cases when $n=2$, and $k=s$ for $n=7$.

For $n=2$, we have $s=t$ and thus by a previous case (where $19\nmid x$ and $n$ even) we see that 
$$(X, Y)=\left(\frac{19^{2(k-t)+1}-1}{2}, \frac{19^{2(k-t)+1}+1}{4}\right).$$
This shows that $(x, y,2)=\left( 19^t\times \frac{19^{2(k-t)+1}-1}{2}, 19^t \times \frac{19^{2(k-t)+1}+1}{4},2 \right)$ is also a solution of \eqref{1} for some non-negative integer $t$.

Again for $n=7$, we have $2s=7t$ and thus $7|s$. So we can write $s=7m$ for some non-negative integer $m$. Therefore by a previous case (where $19\nmid x$ and $n$ odd),
the solution of this equation is
   $(X, Y)=(559, 5)$. This shows that $(x,y,n)=(559\times 19^{7m}, 5\times 19^{2m}, 7)$ is also a solution of \eqref{1} for some integer $m\geq 0$.
   
   Finally
 $2s = \mbox{min} \{2s,2k+1,tn\}$. Then 
 $$
 X^{2}+19^{2(k-s)+1}=4Y^{n}\times 19^{tn-2s}.
 $$
 Considering this equation modulo $19$,  we get $t n = 2 s$. Hence 
 $$
 X^{2}+19^{2(k-s)+1}=4Y^{n}.
 $$
This is same as the last case. Thus we complete the proof.  
\section*{acknowledgement}
The authors are grateful to Prof. Kalyan Chakraborty for his careful reading, helpful comments and suggestions. A. Hoque is supported by the SERB-NPDF (PDF/2017/001958), Govt. of India. R. Sharma  would like to thank to Harish-Chandra Research Institute (HRI) and  Malaviya National Institute of Technology, Jaipur for providing sufficient facility to prepare this manuscript.  The authors would like to thank the referee whose suggestions and comments help to improve the manuscript.

\end{document}